\documentclass[12pt]{article}

\usepackage{amsmath}

\usepackage{natbib}
\usepackage{graphicx}
\usepackage{amssymb,amsmath}
\usepackage{enumerate}
\usepackage{url}

\oddsidemargin
-0.2in \topmargin -.60in
\begin{document}

\title{A Refined Non-asymptotic Tail Bound of Sub-Gaussian Matrix}

\author{
Xianjie~Gao,\quad Hongwei~Zhang\\
School of Mathematical Sciences, Dalian University of Technology \\
Dalian, Liaoning, 116024, P.R. China\\
\texttt{xianjiegao@foxmail.com;hwzhang@dlut.edu.cn}
}
\date{}



\maketitle

\begin{abstract}%
In this paper, we obtain a refined non-asymptotic tail bound for the largest singular value (the soft edge) of sub-Gaussian matrix.  As an application, we use the obtained theorem to compute the tail bound of the Gaussian Toeplitz matrix.

\end{abstract}


\section{Introduction}

Random matrix theory (RMT) has been widely applied in many fields, {\it e.g.,} multivariate statistics [1], high-dimensional data analysis [2], the matrix approximation [3], the combinatorial optimization [4] and the compressed sensing [5]. One main research concern on RMT is to study the tail behavior of the extreme eigenvalues (or singular values) of random matrices.

In general, there are two types of probabilistic statements on the study of probability theory: asymptotic and non-asymptotic. The former aims to analyze the limit behavior of some probability terms, {\it e.g.}, the central limit theorem\begin{equation}\label{eq:clt}
\frac{1}{\sqrt{n}}\sum_{i=1}^nx_i\rightarrow g,\quad (n\rightarrow \infty)
\end{equation}
for Bernoulli random variables $x_1,x_2,\cdots,x_n,\cdots$, where $g$ is Gaussian variable. There have been many well-known asymptotic results on RMT:
\begin{description}
\item[Wigner's semicircle law][6]: Let ${\bf A}_n$ be a $n\times n$ symmetric matrices whose entries are independent Gaussian variables. As dimension $n\rightarrow \infty$,
the spectrum of the Wigner matrices ${\bf W}_n=n^{-1/2}{\bf A}_n$ is distributed according to the semicircle law with density:
\begin{equation}
  f_{sc}(x)=\frac{1}{2\pi}\sqrt{4-x^2}, \quad x\in [-2,2].
\end{equation}

\item[Marchenko-Pastur law][7]: Let ${\bf A}_{m,n}$ $(m\geq n)$ be a $m\times n$ random Gaussian matrix. As the dimensions $m, n\rightarrow \infty$ while the aspect ratio $n/m$ converges to a fix number $y\in (0,1]$, the spectrum of the matrices $\frac{1}{m}{\bf A}^*{\bf A}$ is distributed according to the Marchenko-Pastur law with density:
\begin{equation}\begin{aligned}
 f_{mp}(x)=
  \begin{cases}
\frac{1}{2\pi xy}\sqrt{(b-x)(x-a)}& \text{$a\leq x\leq b$};\\
0& {otherwise},
\end{cases}
\end{aligned}
\end{equation}
where $a=(1-\sqrt{y})^2$ and $b=(1+\sqrt{y})^2$.

\item [Bai-Yin's law][8]: Let ${\bf A}_{m,n}$ $(m\geq n)$ be a $m\times n$ random matrix whose
entries are independent copies of a random variable with zero mean, unit variance, and finite fourth moment. As the dimensions $m, n\rightarrow \infty$ with $n/m$ converging to a fix number $y\in (0,1]$, the $s_{\min}(\bf A)$ and $s_{\max}(\bf A)$ are subjected to Bai-Yin's law:
\begin{align}
  &s_{\min}(\bf A)=\sqrt{\rm m}-\sqrt{\rm n}+{\rm o}(\sqrt{\rm n}),\nonumber\\
  &s_{\max}(\bf A)=\sqrt{\rm m}+\sqrt{\rm n}+{\rm o}(\sqrt{\rm n}),\;\; \mbox{almost surely}.
\end{align}

\end{description}

Although these asymptotic statements can provide a precise limit result when the matrix dimension or sample number goes to the {\it infinity}, they cannot describe in what rate these probability terms converge to their limits. To handle this issue, there arise the non-asymptotic viewpoint to study these probability terms. For example, one of the non-asymptotic statement of the central limit theorem is Hoeffding's inequality:
\begin{equation}\label{eq:ben}
{\mathbb P}\bigg(\frac{1}{\sqrt{n}}\sum_{i=1}^nx_i>t\bigg)\leq 2{\rm e}^{-t^2/2}.
\end{equation}
There have been many research works on RMT from the non-asymptotic viewpoint. Vershynin [9] gave non-asymptotic methods about the properties of sub-Gaussian and sub-exponential matrix.
Tropp [10] proposed a user-friendly framework to study the tail behavior of sums of random matrices. Moreover, there are also other methods for developing the matrix concentration inequalities, {\it e.g.,} exchangeable pairs [11] and Markov chain couplings [12]. To eliminate the dimension dependence of these tail results for random matrices, the intrinsic dimension (or effective dimension) was employed to improve them (see [13],[14]). Recently, Zhang {\it et al.} [15] applied a diagonalization method to obtain the dimension-free tail inequalities of largest singular value for sums of random matrices.

In this paper, we obtain a refined non-asymptotic tail bound for the largest singular value (the soft edge) of sub-Gaussian matrix. We first give a tail bound for the norm of a sub-Gaussian matrix by transforming a sub-Gaussian matrix into a sub-Gaussian variable. We also obtain a tail bound for the norm of a sub-Gaussian matrix by decomposing a sub-Gaussian matrix into a series of sub-Gaussian matrices. By combining the two resulted tail bounds, we obtain the final tail results. As an application, we use the resulted tail inequalities to study the tail behavior of Gaussian Toeplitz matrix.

The rest of this paper is organized as follows. In the next section, we give some preliminary knowledge on random matrices and sub-Gaussian distributions. In Section 3, we present the main results. Section 4 present the application of our results in the study of Gaussian Toeplitz matrix, and the last section concludes paper.

\section{Notations and Preliminaries}\label{sec:pr}

In this section, we give some preliminary knowledge on random matrices and sub-Gaussian distributions.

A random matrix is a matrix whose entries are random variables. Its distribution is characterized by the joint distribution of the entries. The expected value of an $m\times n$ random matrix ${\bf B}$ is the $m\times n$ matrix ${\mathbb E}({\bf B})$ whose entries are the expected values of the corresponding entries of ${\bf B}$, assuming that they all exist.

Let ${\bf B}_{m\times n}$ be a random matrix. Let $S^{n-1}=\{x\in {\mathbb R}^n:\|x\|_2=1\}$ denote the Euclidean sphere in ${\mathbb R}^n$. The largest singular value of ${\bf B}$ is by definition
\begin{equation}\label{eq:op-norm}
  s_{\max}({\bf B})=\|{\bf B}\|=\sup_{x\in {\mathbb R}^n \backslash \{0\}}\frac{\|{\bf B}x\|_2}{\|x\|_2}=\sup_{x\in S^{n-1}}\|{\bf B}x\|_2.
\end{equation}
Given an arbitrary matrix ${\bf B}$, the Hermitian dilation of ${\bf B}$ is defined by
\begin{equation}
  \mathcal{H}(B)=\begin{bmatrix}
0&{\bf B}\\
{\bf B}*& 0
\end{bmatrix}.
\end{equation}
It is ture that $\lambda_{\max}({\mathcal H}{(\bf B)})=\|{\mathcal H}{(\bf B)}\|=\|{\bf B}\|$, where $\lambda_{\max}$ denotes the largest eigenvalue.
The relationship for real function $f$ is the transfer rule. If $f(a)\leq g(a)$ for $a \in I$, then $f({\bf H})\preceq g({\bf H})$ for the eigenvalues of ${\bf H}$ lie in $I$.

Sub-gaussian distributions are referring to a large class of probability
distributions, {\it e.g.,} normal random variables, Bernoulli and all bounded random variables.
\\[8pt]
{\bf Definition 2.1} \ \ {\sl A real-valued random variable $x$ is said to be sub-Gaussian if there exits $c>0$ such that for
every $t>0$
\begin{equation}\label{eq:de-sub-g}
  {\mathbb P}(|x|>t)\leq 2{\rm e}^{-ct^2}.
\end{equation}}

Assuming the sub-Gaussian random variable's mean is zero, the following lemma presents equivalent conditions.
\\[8pt]
{\bf Lemma 2.2} \ \ {\sl Let $x$ be a mean zero (centered) random variable, the following statements are equivalent: 1) $x$ is sub-Gaussian; and 2) $\exists b>0$, $\forall \theta\in {\mathbb R}$, there holds that
    \begin{equation}\label{eq:le-sub-g}
    {\mathbb E}{\rm e}^{\theta x}\leq {\rm e}^{b^2 {\theta}^2/2}.
    \end{equation}}

There are more and more research interests lying in the sub-Gaussian distributions, including spectral properties of random matrices [16] and tail inequalities of sub-Gaussian random vectors [17].

\section{Main Results}\label{sec:psd}
In this section, we obtain a refined upper bound for the largest singular value (the norm) of
sub-Gaussian matrix. We first give a upper bound for the norm of sub-Gaussian matrix by converting into a random sub-Gaussian variable.
\\[8pt]
{\bf Theorem 3.1} \ \ {\sl Let ${\bf B}$ be an $m\times n$ random sub-Gaussian matrix. That is, its entries $x_{ij}$ are i.i.d. centered random variables, obeys the sub-Gaussian distribution. Then there holds that for all $t\geq 0$,
  \begin{equation}
  {\mathbb P}\{\|{\bf B}\|> t\}\leq 2\cdot 5^{(m+n)}\cdot \exp(-ct^2).
 \end{equation}}

The proof of Theorem 3.1 is similar to the Proposition 2.4 of [9], where $m=n$. Here we give the proof of the general case.
\\[8pt]
{\bf Proof}\ \ The main idea of the proof of Theorem 3.1 is to convert the random matrix into a random variable, {\it i.e.,} $\langle{\bf B}x,y\rangle $ is a sub-Gaussian random variable. We then use the covering number to complete the proof.
  \begin{align*}
&{\mathbb P}(\|{\bf B}\|>t)\\
\leq & {\mathbb P}(\max_{\substack{x\in {\mathcal N}\\y\in{\mathcal M}}}\langle {\bf B}x,y\rangle >\frac{t}{4} )\\
\leq & \sum_{\substack{x\in {\mathcal N}\\y\in{\mathcal M}}}{\mathbb P}(\langle{\bf B}x,y\rangle >\frac{t}{4})\\
\leq & |{\mathcal N}||{\mathcal M}|\cdot {\mathbb P}(\langle{\bf B}x,y\rangle >\frac{t}{4})\\
\leq & 2\cdot 5^{(m+n)}\cdot \exp(-ct^2),
\end{align*}
where $\mathcal{N}$, $\mathcal{M}$ are $\frac{1}{2}$-nets of $S^{n-1}$, $S^{m-1}$ respectively, and the bounds on cardinality of the net are $|{\mathcal N}|\leq (1+\epsilon/2)^n$ and $|{\mathcal M}|\leq (1+\epsilon/2)^m$.  \hfill $\blacksquare$

A minor shortcoming of above result is that when the matrix dimension increases, the result becomes very loose. Another method is to obtain the tail bound for matrix sub-Gaussian series. We first introduce the matrix sub-Gaussian moment generating function (mgf) bound.
\\[8pt]
{\bf Proposition 3.2} \ \ {\sl Assume that ${\bf H}$ is a fixed Hermitian matrix and the random variable $x$ obeys the centered sub-Gaussian distribution. Then, there
holds that,
\begin{equation}\label{eq:mgf}
  {\mathbb E}\,{\rm e}^{x\theta {\bf H}}\preceq {\rm e}^{{\theta}^2b^2{\bf H}^2/2}.
\end{equation}}
According to the transfer rule, it is easy to get the proposition. Based on the mgf result (3.2), we develop a tail bound for the matrix sub-Gaussian series.
\\[8pt]
{\bf Theorem 3.3} \ \ {\sl Consider a finite sequence $\{{\bf H}_k: k=1,\ldots,K\}$ of fixed Hermitian matrices with dimension $d$, and $\{x_k:k=1,\ldots,K\}$ be a finite sequence of independent centered sub-Gaussian random variables. Compute the variance parameter
 \begin{equation*}
   \rho:=\|\sum_k{\bf H}_k^2\|.
 \end{equation*}
 Then, for all $t\geq 0$,
 \begin{equation}\label{mids}
 {\mathbb P}\Big\{\lambda_{\max}\Big(\sum_k x_k {\bf H}_k\Big)\geq t\Big\}\leq d\cdot \exp\bigg(-\frac{t^2}{2b^2\rho}\bigg).
\end{equation}}
\\[8pt]
{\bf Proof}\ \ It follows from Proposition 3.2 that, for any $\theta >0$,
\begin{align*}
&{\mathbb P}\Big\{\lambda_{\max}\Big(\sum_k x_k {\bf H}_k\Big)\geq t\Big\}\\
\leq & {\rm e}^{-\theta t}\cdot {\rm tr} \exp\Big(\sum_k \log {\rm E}{\rm e}^{\theta x_k {\bf H}_k}\Big)\\
\leq & {\rm e}^{-\theta t}\cdot {\rm tr} \exp\Big(\frac{{\theta}^2b^2}{2}\sum_k{\bf H}_k^2\Big)\\
\leq & {\rm e}^{-\theta t}\cdot d\cdot \lambda_{\max}\Big( \exp\Big(\frac{{\theta}^2b^2}{2}\sum_k{\bf H}_k^2\Big)\Big)\\
=& d\cdot \exp \Big(-\theta t+\frac{{\theta}^2b^2}{2}\lambda_{\max}\big(\sum_k{\bf H}_k^2\big)\Big)\\
=& d\cdot \exp \Big(-\theta t+\frac{{\theta}^2b^2}{2}\rho\Big),
\end{align*}
where $\rho:=\|\sum_k{\bf H}_k^2\|$, the first inequality follows from Theorem 3.6 of [10]. This inequality holds for any positive $\theta$, so we may take an infimum to complete the
proof. The infimum is attained when $\theta=\frac{t}{b^2\rho}$.  \hfill $\blacksquare$

We apply above result to study the sum of rectangular matrix series by using matrices Hermitian dilation. The following is the general version of Theorem 3.3.
\\[8pt]
{\bf Corollary 3.4} \ \ {\sl Consider a finite sequence $\{{\bf D}_k: k=1,\ldots,K\}$ of fixed matrices with dimension $m\times n$, and $\{x_k:k=1,\ldots,K\}$ be a finite sequence of independent centered sub-Gaussian random variables. Compute the variance parameter
 \begin{equation}
   \rho:=\max \Big\{\Big\|\sum_k{\bf D}_k{\bf D}_k^*\Big\|\, \Big\|\sum_k{\bf D}_k^*{\bf D}_k\Big\|\Big\}.
 \end{equation}
 Then, for all $t\geq 0$,
 \begin{equation}\label{mids}
 {\mathbb P}\Big\{\Big\|\sum_k x_k {\bf D}_k\Big\|\geq t\Big\}\leq (m+n)\cdot \exp\bigg(-\frac{t^2}{2b^2\rho}\bigg).
\end{equation}}
\\[8pt]
{\bf Proof}\ \ According to Hermitian dilation we know that
  \begin{equation*}
  \Big\|\sum_k x_k {\bf D}_k\Big\|
  =\lambda_{\max}\Big({\mathcal H}{\Big(\sum_k x_k {\bf D}_k\Big)}\Big)
  =\lambda_{\max}\Big(\sum_k x_k{\mathcal H}({{\bf D}_k})\Big).
  \end{equation*}

We invoke Theorem 3.3 to obtain the tail bound for the sum of rectangular matrix series. The matrix variance parameter $\rho$ satisfies the relation:
\begin{equation*}
 \rho=\Big\|\sum_k{\mathcal H}({{\bf D}_k)^2}\Big\|
 =\left\|\begin{matrix}
     \sum_k{\bf D}_k{\bf D}_k^* & 0 \\
     0 & \sum_k{\bf D}_k^*{\bf D}_k
\end{matrix}\right\|
 =\max \Big\{\Big\|\sum_k{\bf D}_k{\bf D}_k^*\Big\|\, \Big\|\sum_k{\bf D}_k^*{\bf D}_k\Big\|\Big\}
\end{equation*}
This completes the proof.  \hfill $\blacksquare$

Based on the general version of tail bound for matrix sub-Gaussian series, we obtain another tail bound for the norm of the sub-Gaussian matrix.
\\[8pt]
{\bf Theorem 3.5} \ \ {\sl Under the notations and conditions in Theorem 3.1. Then there holds that for all $t\geq 0$,
  \begin{equation}\label{eq:matrix-sub-g}
  {\mathbb P}\{\|{\bf B}\|> t\}\leq (m+n)\cdot \exp\Big(-\frac{t^2}{2b^2m}\Big).
\end{equation}}
\\[8pt]
{\bf Proof}\ \ In order to use Corollary 3.4, we decompose matrix as a matrix sub-Gaussian series:
 \begin{equation*}
   {\bf B}=\sum_{ij}x_{ij}{\bf E}_{ij}, \quad i=1,\ldots,m, \quad j=1,\ldots,n.
 \end{equation*}
 The matrix ${\bf E}_{ij}$ has a element one in the $(i,j)$ position and zeros elsewhere.
By calculating $\rho=m$, the conclusion is established by using Corollary 3.4.   \hfill $\blacksquare$


The combination of Theorem 3.1 and Theorem 3.5 leads to the following refined upper bound for the largest singular value (the soft edge) of sub-Gaussian matrix.
\\[8pt]
{\bf Theorem 3.6} \ \ {\sl Follow the notations and conditions in Theorem 3.1. Then there holds that for all $t\geq 0$,
 \begin{equation} \begin{aligned}\label{eq:refine}
  {\mathbb P}\{\|{\bf B}\|> t\}\leq
  \begin{cases}
(m+n)\cdot \exp\Big(-\frac{t^2}{2b^2m}\Big)& \text{$0<t\leq \sqrt{\frac{2b^2m}{1-2b^2mc}\log \frac{m+n}{2\cdot 5^{m+n}}}$};\\
2\cdot 5^{(m+n)}\cdot \exp(-ct^2)& \text{$t> \sqrt{\frac{2b^2m}{1-2b^2mc}\log \frac{m+n}{2\cdot 5^{m+n}}}$}.
\end{cases}
\end{aligned}
\end{equation}}

\baselineskip 15pt
\section{Application: Gaussian Toeplitz Matrix}
In this section, we use our theoretical findings to compute the tail bound of the Gaussian Toeplitz matrix. The Gaussian Toeplitz matrix is an example of Gaussian random matrix which has been widely used in various fields, {\it e.g.,}differential equations, spline functions, and signal processing [18]. We consider a unsymmetric Gaussian Toeplitz matrix ${\bf T}\in\mathbb{C}^{d\times d}$ in the following form:
\begin{equation}
{\bf T}=\begin{bmatrix}
\gamma_{0}&\gamma_{1}&\gamma_{2} &\cdots&\gamma_{d-1}\\
\gamma_{-1}&\gamma_{0}&\gamma_{1}&\cdots &\gamma_{d-2} \\
\gamma_{-2}&\gamma_{-1}&\gamma_{0}&\cdots &\gamma_{d-3}\\
\vdots&\vdots&\vdots &\ddots&\vdots\\
\gamma_{-(d-1)}&\gamma_{-(d-2)}&\gamma_{-(d-3)}&\cdots&\gamma_{0}
\end{bmatrix},
\end{equation}
where $\gamma_{-(d-1)},\ldots,\gamma_{d-1}$ are independent standard normal variables. The Gaussian Toeplitz matrix ${\bf T}$ can be represented as a matrix Gaussian series:
\begin{equation}
 {\bf T}=\gamma_0{\bf I}+\sum_{j=1}^{d-1}\gamma_{j}{\bf C}^j+\sum_{j=1}^{d-1}\gamma_{-j}({\bf C}^j)^T,
\end{equation}
where ${\bf C}^j$ is the $j$-th power of ${\bf C}$ with
\begin{equation*}
{\bf C}=\begin{bmatrix}
0&1& & &\\
&0&1& & \\
&&\ddots&\ddots &\\
&&&0&1\\
&&&&0
\end{bmatrix}.
\end{equation*}

Using the Theorem 3.6, we can compute the tail bound of the Gaussian Toeplitz matrix. First, we calculate
\begin{align*}
  ({\bf C}^j)({\bf C}^j)^T=\sum_{k=1}^{d-j}{\bf E}_{kk} \quad and \quad ({\bf C}^j)^T({\bf C}^j)=\sum_{k=j+1}^{d}{\bf E}_{kk}.
\end{align*}
The matrix variance parameter $\rho=d$ can be calculated:
\begin{align*}
  {\bf I}^2+\sum_{j=1}^{d-1}({\bf C}^j)({\bf C}^j)^T+\sum_{j=1}^{d-1}({\bf C}^j)^T({\bf C}^j)=d{\bf I}_d.
\end{align*}
For Gaussian matrix, $b=1$, $c=\frac{1}{2}$. Through the application of Theorem 3.6, tail bound of the Gaussian Toeplitz matrix is presented, for all $t\geq 0$,
\begin{align}\label{eq:refine}
  {\mathbb P}\{\|{\bf T}\|> t\}\leq
  \begin{cases}
2d\cdot \exp\Big(-\frac{t^2}{2d}\Big)& \text{$0<t\leq \sqrt{\frac{2d}{1-2d}\log \frac{2d}{2\cdot 5^{d}}}$};\\
2\cdot 5^{d}\cdot \exp(-\frac{t^2}{2})& \text{$t> \sqrt{\frac{2d}{1-2d}\log \frac{2d}{2\cdot 5^{d}}}$}.
\end{cases}
\end{align}

\section{Conclusion}

In this paper, we first present the tail bounds for the largest singular value of sub-Gaussian matrix and matrix sub-Gaussian series. We then obtain a refined non-asymptotic tail bound for the largest singular value (the soft edge) of sub-Gaussian matrix. As an application, we finally compute the tail bound of Gaussian Toeplitz matrix.


\end{document}